\documentclass{birkjour}
 \newtheorem{thm}{Theorem}[section]

 \theoremstyle{definition}
 
 \theoremstyle{remark}
 
 \newtheorem{ex}{Example}[section]
 \numberwithin{equation}{section}
\begin{document}

\title[Ricci solitons and certain related metrics]
{Ricci solitons and certain related metrics on 3-dimensional trans-Sasakian manifold}

\author[S. Sarkar]{Sumanjit Sarkar}
\address{Department of Mathematics\\
Jadavpur University\\
Kolkata-700032, India}
\email{imsumanjit@gmail.com}

\thanks{The first author is the corresponding author and this work was financially supported by UGC Senior Research Fellowship of India, Sr. No. 2061540940. Ref. No:21/06/2015(i)EU-V.}

\author[S. Dey]{Santu Dey}
\address{Department of Mathematics\\
Bidhan Chandra College\\
Asansol, Burdwan, West Bengal-713304, India}
\email{santu.mathju@gmail.com}

\author[A. Bhattacharyya]{Arindam Bhattacharyya}
\address{Department of Mathematics\\
Jadavpur University\\
Kolkata-700032, India}
\email{bhattachar1968@yahoo.co.in}

\subjclass{53C25, 53C44, 53D15}
\keywords{Ricci flow, Ricci soliton, $*$-conformal Ricci soliton, trans-Sasakian manifold.}

\begin{abstract}In this paper we study certain types of metrics such as Ricci soliton, $*$-conformal Ricci soliton in 3-dimensional trans-Sasakian manifold. First we have shown that a 3-dimensional trans-Sasakian manifold of type $(\alpha,\beta)$ admits a Ricci soliton where the covariant derivative of potential vector field in the direction of unit vector field $\xi$ is orthogonal to $\xi$. It is also shown that if the structure functions satisfy $\alpha^2=\beta^2$ then the covariant derivative of the potential vector field in the direction of $\xi$ is a constant multiple of $\xi$. Further, we have evolved the nature of scalar curvature when the manifold satisfies $*$-conformal Ricci soliton of type $(\alpha,\beta)$, provided $\alpha \neq 0$. Finally, we present an example to verify our findings.
\end{abstract}
\maketitle

\section{\textbf{Introduction}}
\medskip
Richard S. Hamilton  introduced the concept of Ricci flow (for details see \cite{Ham}) which was named after great Italian mathematician Gregorio Ricci-Curbastro. Later Grigori Perelman found it very useful to solve Poincare conjecture. If we take a smooth closed (compact without boundary) Riemannian manifold $M$ equipped with a smooth Riemannian metric $g$ then the Ricci flow is defined by the geometric evolution equation,
\begin{equation}\label{1.1}
  \frac{\partial g(t)}{\partial t}=-2S(g(t))
\end{equation}
where $S$ is the Ricci curvature tensor of the manifold and $g(t)$ is a one-parameter family of metrices on M.\par
A Riemannian manifold $(M,g)$ is called a Ricci soliton if there exists a vector field $V$ and a constant $\lambda$ such that the following equation holds,
\begin{equation}\label{1.2}
  \frac{1}{2}\mathcal{L}_Vg+S+\lambda g=0
\end{equation}
where $\mathcal{L}_V$ denotes Lie derivative along the direction $V$ and $\lambda$ is a non-zero constant. The vector field $V$ is called potential vector field and $\lambda$ is called soliton constant. Ricci soliton which is a natural extension of Einstein manifold is a self-similar solution of Ricci flow. The potential vector field $V$ and soliton constant $\lambda$ play  vital roles while determining the nature of the soliton. A soliton is said to be shrinking, steady or expanding according as $\lambda<0$, $\lambda=0$ or $\lambda>0$. Now if $V$ is Killing then the Ricci soliton reduces to Einstein manifold. Compact Ricci solitons are the fixed points of the Ricci flow (\ref{1.1}) projected from the space of metrics onto its quotient modulo diffeomorphisms and scalings, and often arise as blow-up limits for the Ricci flow on compact manifolds.\par
In 2005, A. E. Fischer \cite{Fischer} has introduced conformal Ricci flow which is a mere generalisation of the classical Ricci flow equation that modifies the unit volume constraint to a scalar curvature constraint. The conformal Ricci flow equation was given by,
\begin{eqnarray*}
  \frac{\partial g}{\partial t}+2(S+\frac{g}{n}) &=& -pg, \\
  r(g) &=& -1,
\end{eqnarray*}
where $r(g)$ is the scalar curvature of the manifold, $p$ is scalar non-dynamical field and $n$ is the dimension of the manifold. Corresponding to the conformal Ricci flow equation in 2015, N. Basu and A. Bhattacharyya \cite{Basu} introduced the notion of conformal Ricci soliton equation as a generalization of Ricci soliton equation given by,
\begin{equation}\label{1.3}
\mathcal{L}_Vg+2S+[2\lambda-(p+\frac{2}{n})]g=0.
\end{equation}
In 2014, Kaimakamis and Panagiotidou \cite{Kai} modified the definition of Ricci soliton where they have used $*$-Ricci tensor $S^*$ which was introduced by Tachibana \cite{Tachi} and Hamada \cite{hama} respectively, in place of Ricci tensor $S$. The $*$-Ricci tensor $S^*$ is defined by
\begin{equation*}
  S^*(X,Y)=\frac{1}{2}(trace\{\phi.R(X,\phi Y)\})
\end{equation*}
for all vector fields $X$ and $Y$ on $M$. They have used the concept of $*$-Ricci soliton within the framework of real hypersurfaces of a complex space form. A pseudo-Riemannian metric $g$ is called a $*$-Ricci soliton if there exists a constant $\lambda$ and a vector field $V$ such that
\begin{equation*}
  \mathcal{L}_Vg+2S^*+2\lambda g=0.
\end{equation*}
Note that, $*$-soliton is trivial if the vector field $V$ is Killing, and in this case, the manifold becomes $*$-Einstein. Thus, it is considered as a natural generalization of $*$-Einstein metric. A $*$-Ricci soliton is said to be almost $*$-Ricci soliton if $\lambda$ is a smooth function on $M$. Moreover, an almost $*$-Ricci soliton is called shrinking, steady and expanding according to as $\lambda$ is positive, zero and negative, respectively. In \cite{Kai}, it was studied a real hypersurfaces of a non-flat complex space form admitting a $*$-Ricci soliton whose potential vector field is the structure vector field and proved that a real hypersurface in a complex projective space does not admit a $*$-Ricci soliton. They have also shown that a real hypersurface of complex hyperbolic space admitting a $*$-Ricci soltion is locally congruent to a geodesic hypersphere.\\

Further P. Majhi and D. Dey \cite{Majhi} in 2020 revised the aforementioned definition of $*$-Ricci soliton with the help of (\ref{1.3}) and defined $*$-Conformal Ricci soliton as
\begin{equation}\label{1.4}
  \mathcal{L}_Vg+2S^*+[2\lambda-(p+\frac{2}{n})]g=0.
\end{equation}
Ricci solitons have been studied in many contexts: on K$\ddot{a}$hler manifolds\cite{Cho}, on contact and Lorentzian manifolds \cite{Bag}, \cite{Ash}, on K-contact manifolds \cite{Sharma} etc. by many authors. Later, Nagaraja and Premalatha \cite{Nag} studied the nature of Ricci soliton on Kenmotsu manifold; C$\check{a}$lin and Crasmareanu \cite{Cal} on f-Kenmotsu manifold; He and Zhu \cite{He} on Sasakian manifold and Ingalahalli and Bagewadi \cite{Ing} on $\alpha$-Sasakian manifold. Recently, in 2017, Yaning Wang \cite{Wang} proved that if a 3-dimensional cosymplectic manifold $M^3$ admits a Ricci soliton, then either $M^3$ is locally flat or the potential vector field is an infinitesimal contact transformation. Also, S. Pahan and A. Bhattacharyya gave some insight on trans-Sasakian manifold in \cite{Sampa}. In 2016, T. Dutta et al. studied conformal Ricci soliton on 3-dimensional trans-Sasakian manifold\cite{Dutta}.\\

As follows in the literature, $*$-Ricci soliton on contact geometry studied by many authors: on Sasakian and $(\kappa,\mu)$-contact manifold by Ghosh and Patra \cite{Ghosh-Patra}, on $(\kappa,\mu)^\prime$-almost Kenmotsu manifolds by Dai, Zhao and De\cite{Dai}, on contact 3-manifolds by Wang \cite{WangY} etc. It is worthy to mention that in \cite{Dey} Dey and Majhi considered $*$-Ricci soliton on 3-dimensional trans-Sasakian manifold and proved that if the metric of the manifold represents $*$-Ricci soliton and if it satisfies a certain condition then the manifold reduces to a $\beta$-Kenmotsu manifold. Also, very recently generalizations of $*$-Ricci soliton on contact geometry have been studied by \cite{dey1, dey2}.

Motivated from above mentioned well praised works we have studied behaviour of Ricci soliton and $*$-conformal Ricci soliton on 3-dimensional trans-Sasakian manifold. In the later sections we have revisited some definitions and important properties of 3-dimensional trans-Sasakian manifold and after that the main result of this paper, containing Theorem 3.1 and Theorem 3.2 has been described. We also provide an example to justify our findings.
\medskip

\section{\textbf{Preliminaries}}
By \cite{Blair}, a differentiable manifold $M$ of dimension $2n+1$ is said to be an almost contact structure or $(\phi,\xi,\eta)$ structure if $M$ admits a $(1,1)$ tensor field $\phi$, a vector field $\xi$, a 1-form $\eta$ satisfying:
\begin{eqnarray}
  \phi ^2 &=& -I + \eta \otimes \xi \label{2.1}\\
  \eta (\xi ) &=& 1\label{2.2}
\end{eqnarray}
where $I$ is the identity mapping.
A Riemannian metric $g$ is said to be \emph{compatible metric} if it satisfies:
\begin{equation}\label{2.3}
  g( \phi X , \phi Y) = g(X,Y)- \eta (X) \eta (Y)
\end{equation}
A manifold having almost contact structure along with compatible Riemannian metric is called \emph{almost contact metric manifold}.\par
In an almost contact metric manifold the following conditions are satisfied:
\begin{eqnarray}
  \phi \xi &=& 0\label{2.4} \\
  \eta \circ \phi &=& 0 \label{2.5}\\
  g(X, \xi ) &=& \eta (X) \label{2.6}\\
  g(\phi X,Y) &=& -g(X,\phi Y)\label{2.7}
\end{eqnarray}
for arbitrary $X,Y \in \chi (M)$.\par
Let $M$ be a $2n+1$ dimensional almost contact manifold. Then we define an almost complex structure $J$ on $M\times\mathbb{R}$ by $J(X,f\frac{d}{dt})=(\phi X-f\xi ,\eta (X)\frac{d}{dt})$, where $X$ is a tangent to $M$, $t$ is the coordinate on $\mathbb{R}$ and $f$ a $C^\infty$ function on $M\times\mathbb{R}$. Clearly $J^2=-I$. If $J$ is integrable then the almost contact structure is said to be normal. The normality of an almost contact metric manifold is equivalent with vanishing of the tensor field $[\phi,\phi]+2d\eta\otimes\xi$, where $[\phi,\phi]$ is the Nijenhuis torsion tensor of $\phi$(for more details see \cite{Blair}).\par
An almost contact metric manifold $M$ is called a trans-Sasakian manifold if $(M\times\mathbb{R},J,G)$, where $G$ is the product metric on $M\times\mathbb{R}$, belongs to the class $W_4$ (see \cite{Gray}). If there are smooth functions $\alpha,\beta$ on an almost contact metric manifold $(M,\phi,\xi,\eta,g)$ satisfying:
\begin{eqnarray}
  (\nabla_X \phi)Y&=&\alpha[g(X,Y)\xi-\eta(Y)X]+\beta[g(\phi X,Y)\xi-\eta(Y)\phi X]\label{2.8}\\
  (\nabla_X\eta)(Y)&=&-\alpha g(\phi X,Y)+\beta g(\phi X,\phi Y)\label{2.9}
\end{eqnarray}
where $X,Y\in\chi(M)$ are arbitrary and $\nabla$ is the Levi-Civita connection of $g$ on $M$, then the manifold is called trans-Sasakian manifold of type $(\alpha,\beta)$. $\alpha, \beta$ are called structure functions of the manifold. Trans-Sasakian manifolds of type $(0,0),(\alpha,0),(0,\beta)$ are called cosymplectic, $\alpha$-Sasakian, $\beta$-Kenmotsu manifolds respectively. Then form (\ref{2.9}) we can deduce that:
\begin{equation}\label{2.10}
  \nabla_X\xi=-\alpha\phi X+\beta(X-\eta(X)\xi)
\end{equation}\par
 Marrero \cite{Marr} showed that a trans-Sasakian manifold of dimension $\geq 5$ is either cosymplectic or $\alpha$-Sasakian or $\beta$-Kenmotsu. So proper trans-Sasakian manifold exists for dimension 3. In a 3 dimensional trans-Sasakian manifold the following relations hold:\par
 \begin{eqnarray}\label{2.11}
 R(X,Y)Z &=& (\frac{r}{2} + 2\xi\beta - 2(\alpha^2 - \beta^2))(g(Y,Z)X - g(X,Z)Y)-g(Y,Z) \nonumber\\
&&[(\frac{r}{2}+\xi\beta-3(\alpha^2-\beta^2))\eta(X)\xi-\eta(X)(\phi D\alpha-D\beta)+(X\beta\nonumber\\
  &&+(\phi X)\alpha)\xi]+g(X,Z)[(\frac{r}{2}+\xi\beta-3(\alpha^2-\beta^2))\eta(Y)\xi-\eta(Y)\nonumber\\
  &&(\phi D\alpha-D\beta)+(Y\beta+(\phi Y)\alpha)\xi)]- [(Z\beta+(\phi Z)\alpha)\eta(Y)+\nonumber\\
 &&(Y\beta+(\phi Y)\alpha)\eta(Z)+(\frac{r}{2}+\xi\beta-3(\alpha^2-\beta^2))\eta(Y)\eta(Z)]X\nonumber\\
  &&+[(Z\beta+(\phi Z)\alpha)\eta(X)+(X\beta+(\phi X)\alpha)\eta(Z)+(\frac{r}{2}+\xi\beta-\nonumber\\
  &&3(\alpha^2-\beta^2))\eta(X)\eta(Z)]Y
 \end{eqnarray}\\
where $Df$ denotes the gradient of the smooth function $f$ defined on $M$.
\begin{eqnarray}
S(X,Y)&=&(\frac{r}{2}+\xi\beta-(\alpha^2-\beta^2))g(X,Y)-(\frac{r}{2}+\xi\beta-3(\alpha^2-\beta^2))\eta(X)\nonumber\\
&&\eta(Y)-(Y\beta+(\phi Y)\alpha)\eta(X)-(X\beta+(\phi X)\alpha)\eta(Y)\label{2.12}\\
S(X,\xi)&=&(2(\alpha^2-\beta^2)-\xi\beta)\eta(X)-X\beta-(\phi X)\alpha\label{2.13}
\end{eqnarray}
where $R,S,r$ are Riemannian curvature tensor, Ricci tensor of type (0,2) and scalar curvature of the manifold respectively and $\alpha,\beta$ are smooth functions on the manifold (for details see \cite{Sampa}).\par
Here in this paper we restricted the smooth functions $\alpha,\beta$ to be constant functions. Then we got some special relations compatible to our restrictions:
\begin{eqnarray}
  R(X,Y)\xi &=& (\alpha^2-\beta^2)(\eta(Y)X-\eta(X)Y) \label{2.14}\\
  \nonumber S(X,Y) &=& (\frac{r}{2}-(\alpha^2-\beta^2))g(X,Y)-(\frac{r}{2}-3(\alpha^2-\beta^2))\eta(X)\eta(Y)\\
  &&\label{2.15}\\
  S(X,\xi) &=& 2(\alpha^2-\beta^2)\eta(X)\label{2.16}\\
    QX &=& (\frac{r}{2}-(\alpha^2-\beta^2))X-(\frac{r}{2}-3(\alpha^2-\beta^2))\eta(X)\xi\label{2.17}
\end{eqnarray}
where $Q$ is the Ricci operator given by $S(X,Y)=g(QX,Y)$. The expression of $*$-Ricci tensor (for details see Lemma 3.1 of \cite{Dey}) on a 3-dimensional trans-Sasakian manifold is given by

\begin{equation}\label{2.18}
  S^*(X,Y)=\frac{1}{2}(r-4(\alpha^2-\beta^2))[g(X,Y)-\eta(X)\eta(Y)]
\end{equation}

\section{\textbf{Main Results}}
\medskip
In this section we consider the metric of 3-dimensional trans-Sasakian manifold as a Ricci soliton and a $*$-conformal Ricci soliton and proved the following two results.
\begin{thm}
Let $M$ be a 3-dimensional trans-Sasakian manifold of type $(\alpha,\beta)$ admitting a Ricci soliton where the structure functions $\alpha$ and $\beta$ are constant. Then the following relations are satisfied:\\
(i) If $\nabla_\xi V$ is orthogonal to $\xi$, then the soliton is shrinking for $\alpha^2<\beta^2$, steady for $\alpha^2=\beta^2$ and expanding for $\alpha^2>\beta^2$.\\
(ii) If $\alpha^2=\beta^2$, then the covariant derivative of the potential vector field $V$ in the direction of $\xi$ is a constant multiple of $\xi$.
  \medskip
  \begin{proof}
    In a 3-dimensional trans-Sasakian manifold where $\alpha$ and $\beta$ are non-zero constant, we know from (\ref{2.17}) that the Ricci operator can be written as,
    \begin{equation}\label{3.1}
      QX=(\frac{r}{2}-(\alpha^2-\beta^2))X-(\frac{r}{2}-3(\alpha^2-\beta^2))\eta(X)\xi
    \end{equation}
    Where $r$ denotes the scalar curvature of the manifold and $X\in \chi(M)$ is any vector field. The aforementioned equation implies that it is an $\eta$-Einstein manifold. Now taking covariant derivative of (3.1) w.r.t arbitrary $Y\in \chi(M)$ we have,
    \begin{eqnarray}\label{3.2}
      (\nabla_YQ)X &=& \frac{1}{2}(Yr)X-\frac{1}{2}(Yr)\eta(X)\xi-(\frac{r}{2}-3(\alpha^2-\beta^2))[-\alpha g(\phi Y,X)\xi\nonumber\\
      && +\beta g(X,Y)\xi-\alpha\eta(X)\phi Y+\beta\eta(X)Y-2\beta\eta(X)\eta(Y)\xi]
    \end{eqnarray}
      Contracting $X$ and using the well-known formula $trace\{X\rightarrow(\nabla_XQ)Y\}=\frac{1}{2}(Yr)$ in (\ref{3.2}) we have,
      \begin{equation}\label{3.3}
        \xi r=-2r\beta+12(\alpha^2-\beta^2)\beta
      \end{equation}
  Using (2.15) in the definition of Ricci soliton (\ref{1.2}) we have for any vector fields $Y,Z\in\chi(M)$,
  \begin{equation}\label{3.4}
    (\mathcal{L}_Vg)(Y,Z)=(2\lambda-r+2(\alpha^2-\beta^2))g(Y,Z)+(r-6(\alpha^2-\beta^2))\eta(Y)\eta(Z)
  \end{equation}
  Now taking covariant derivative of (\ref{3.4}) along an arbitrary vector field $X\in\chi(M)$,
\begin{eqnarray}\label{3.5}
  (\nabla_X\mathcal{L}_Vg)(Y,Z) &=& -(Xr)g(Y,Z)+(Xr)\eta(Y)\eta(Z)+(r-6(\alpha^2-\beta^2)) \nonumber\\
  &&\{-\alpha g(\phi X,Y)\eta(Z)+\beta g(\phi X,\phi Y)\eta(Z)-\alpha g(\phi X,Z)\nonumber\\
  && \eta(Y)+\beta g(\phi X,\phi Z)\eta(Y)\}
\end{eqnarray}
Again for any vector fields $X,Y,Z\in\chi(M)$ we know (for further information please refer to \cite{Yano}),
\begin{eqnarray}
  (\mathcal{L}_V\nabla_Xg-\nabla_X\mathcal{L}_Vg-\nabla_{[V,X]}g)(Y,Z)&=&-g((\mathcal{L}_V\nabla)(X,Y),Z)\nonumber\\
  &&-g((\mathcal{L}_V\nabla)(X,Z),Y)\nonumber
\end{eqnarray}
 Since $g$ is Riemannian metric connection, $\nabla g=0$. So the above equation reduces to,
\begin{eqnarray*}
  (\nabla_X\mathcal{L}_Vg)(Y,Z)=g((\mathcal{L}_V\nabla)(X,Y),Z)+g((\mathcal{L}_V\nabla)(X,Z),Y)
\end{eqnarray*}
 Using symmetry of $\mathcal{L}_V\nabla$ i.e., $(\mathcal{L}_V\nabla)(X,Y)=(\mathcal{L}_V\nabla)(Y,X)$ we can have,
\begin{eqnarray*}
   2g((\mathcal{L}_V\nabla)(X,Y),Z)=(\nabla_X\mathcal{L}_Vg)(Y,Z)+(\nabla_Y\mathcal{L}_Vg)(Z,X)-(\nabla_Z\mathcal{L}_Vg)(X,Y)
 \end{eqnarray*}
 Using (\ref{3.5}) in the above equation we get,
 \begin{eqnarray}\label{3.6}
   (\mathcal{L}_V\nabla)(X,Y) &=& -\frac{1}{2}(Xr)Y-\frac{1}{2}(Yr)X+\frac{1}{2}g(\phi X,\phi Y)Dr+\frac{1}{2}(Xr)\eta(Y)\xi\nonumber \\
    && +\frac{1}{2}(Yr)\eta(X)\xi+(r-6(\alpha^2-\beta^2))\{-\alpha\eta(Y)\phi X-\alpha\eta(X) \nonumber \\
    && \phi Y+\beta g(\phi X,\phi Y)\xi\}
 \end{eqnarray}
 Taking covariant derivative w.r.t arbitrary vector field we have,
 \begin{eqnarray}
   (\nabla_X\mathcal{L}_V\nabla)(Y,Z) &=& -\frac{1}{2}g(Z,\nabla_XDr)Y-\frac{1}{2}g(Y,\nabla_XDr)Z+\frac{1}{2}g(\phi Y,\phi Z)\nonumber\\
   &&(\nabla_XDr)-\alpha\eta(Z)(Xr)\phi Y-\alpha\eta(Y)(Xr)\phi Z+\frac{1}{2}\{(Zr)\nonumber\\
   &&\eta(Y)+(Yr)\eta(Z)\}(\nabla_X\xi)+\frac{1}{2}\{g(Y,\nabla_XDr)\eta(Z)-\alpha\nonumber\\
   &&(Yr)g(\phi X,Z)+\beta g(\phi X,\phi Z)(Yr)+g(Z,\nabla_XDr)\eta(Y)\nonumber\\
   &&-\alpha(Zr)g(\phi X,Y)+\beta g(\phi X,\phi Y)(Zr)+2\beta g(\phi Y,\phi Z)\nonumber\\
   &&(Xr)\}\xi+\frac{1}{2}\{\alpha g(\phi X,Y)\eta(Z)-\beta g(\phi X,\phi Y)\eta(Z)+\alpha\nonumber\\
   && g(\phi X,Z)\eta(Y)-\beta g(\phi X,\phi Z)\eta(Y)\}Dr+(r-6(\alpha^2-\nonumber\\
   &&\beta^2))[\{\alpha^2g(\phi X,Z)-\alpha\beta g(\phi X,\phi Z)\}\phi Y+\{\alpha^2g(\phi X,Y)\nonumber\\
   &&-\alpha\beta g(\phi X,\phi Y)\}\phi Z-\alpha\eta(Z)((\nabla_X\phi)Y)-\alpha\eta(Y)\nonumber\\
   &&((\nabla_X\phi)Z)+\beta g(\phi Y,\phi Z)(\nabla_X\xi)+\{\alpha\beta g(\phi X,Y)\eta(Z)\nonumber\\
   &&-\beta^2g(\phi X,\phi Y)\eta(Z)+\alpha\beta g(\phi X,Z)\eta(Y)\nonumber\\
   &&-\beta^2g(\phi X,\phi Z)\eta(Y)\}\xi]\nonumber
   \end{eqnarray}
From Yano (for more details see \cite{Yano}) we know $(\mathcal{L}_VR)(X,Y)Z=(\nabla_X\mathcal{L}_V\nabla)\\(Y,Z)-(\nabla_Y\mathcal{L}_V\nabla)(X,Z)$. Using this formula in the above equation we get,
\begin{eqnarray}\label{3.7}
  (\mathcal{L}_VR)(X,Y)Z &=& \frac{1}{2}g(Z,\nabla_YDr)X-\frac{1}{2}g(Z,\nabla_XDr)Y+\frac{1}{2}\{-\alpha(Yr)g(\phi X,\nonumber\\
  &&Z)-\beta(Yr)g(\phi X,\phi Z)+g(Z,\nabla_XDr)\eta(Y)-\alpha(Zr)\nonumber\\
  &&g(\phi X,Y)+\alpha(Xr)g(\phi Y,Z)+\beta g(\phi Y,\phi Z)(Xr)-g(Z,\nonumber\\
  &&\nabla_YDr)\eta(X)+\alpha(Zr) g(\phi Y,X)\}\xi+\alpha\eta(Z)(Yr)\phi X-\nonumber\\
  &&\alpha\eta(Z)(Xr)\phi Y+\alpha\{\eta(X)(Yr)-\eta(Y)(Xr)\}\phi Z+ \nonumber\\
  &&\frac{1}{2}\{\alpha g(\phi X,Y)\eta(Z)+\alpha g(X,\phi Y)\eta(Z)-\alpha g(\phi X,Z) \nonumber\\
  &&\eta(Y)-\beta g(\phi X,\phi Z)\eta(Y)-\alpha g(\phi Y,Z)\eta(X)+\beta g(\phi Y,\nonumber\\
   &&\phi Z)\eta(X)\}Dr+\frac{1}{2}\{(Yr)\eta(Z)+(Zr)\eta(Y)\}(\nabla_X\xi)-\nonumber\\
  &&\frac{1}{2}\{(Xr)\eta(Z)+(Zr)\eta(X)\}(\nabla_Y\xi)+\frac{1}{2}g(\phi Y,\phi Z)(\nabla_XDr)\nonumber\\
  &&-\frac{1}{2}g(\phi X,\phi Z)(\nabla_YDr)+(r-6(\alpha^2-\beta^2))[\{\alpha\beta g(\phi Y,\phi Z)\nonumber\\
  &&-\alpha^2g(\phi Y,Z)\}\phi X-\{\alpha\beta g(\phi X,\phi Z)-\alpha^2g(\phi X,Z)\}\phi Y+\nonumber\\
  &&2\alpha^2g(\phi X,Y)\phi Z+\{2\alpha\beta g(\phi X,Y)\eta(Z)+\alpha\beta g(\phi X,Z) \nonumber\\
   &&\eta(Y)-\beta^2g(\phi X,\phi Z)\eta(Y)-\alpha\beta g(\phi Y,Z)\eta(X)+\beta^2\nonumber\\
  &&g(\phi Y,\phi Z)\eta(X)\}\xi+\beta g(\phi Y,\phi Z)(\nabla_X\xi)-\beta g(\phi X,\phi Z)\nonumber\\
  &&(\nabla_Y\xi)-\alpha\eta(Z)((\nabla_X\phi)Y)-\alpha\eta(Y)((\nabla_X\phi)Z)+\alpha\eta(Z)\nonumber\\
  &&((\nabla_Y\phi)X)+\alpha\eta(X)((\nabla_Y\phi)Z)
\end{eqnarray}
This equation holds for any $X,Y,Z\in\chi(M)$. Contracting $X$ in (\ref{3.7}) we get,
\begin{equation}\label{3.8}
  (\mathcal{L}_VS)(Y,Z)=(\frac{\Delta r}{2}-6\alpha^4+12\alpha^2\beta^2-6\beta^4+r\alpha^2-r\beta^2)g(\phi Y,\phi Z)
\end{equation}
for any $Y,Z\in\chi(M)$. Again from (\ref{2.15}) we can have,
\begin{eqnarray}\label{3.9}
  (\mathcal{L}_VS)(Y,Z)&=& \frac{1}{2}g(\phi Y,\phi Z)(Vr)+(\frac{r}{2}-(\alpha^2-\beta^2))\{g(\nabla_YV,Z)\nonumber\\
  &&+g(Y,\nabla_ZV)\}-(\frac{r}{2}-3(\alpha^2-\beta^2))\{\eta(Z)((\nabla_V\eta)Y)\nonumber\\
  &&+\eta(Y)((\nabla_V\eta)Z)+\eta(Z)\eta(\nabla_YV)+\eta(Y)\eta(\nabla_ZV)\}\nonumber\\
  &&
\end{eqnarray}
Comparing (\ref{3.8}) with (\ref{3.9}) yields,
\begin{eqnarray}\label{3.10}
  &&(\frac{\Delta r}{2}-6\alpha^4+12\alpha^2\beta^2-6\beta^4+r\alpha^2-r\beta^2)g(\phi Y,\phi Z)=\frac{1}{2}\{g(\phi Y,\phi Z)(Vr)\nonumber\\
  &&+(\frac{r}{2}-(\alpha^2-\beta^2))\{g(\nabla_YV,Z)+g(Y,\nabla_ZV)\}-(\frac{r}{2}-3(\alpha^2-\beta^2))\{\eta(Z)\nonumber\\
  &&((\nabla_V\eta)Y)+\eta(Y)((\nabla_V\eta)Z)+\eta(Z)\eta(\nabla_YV)+\eta(Y)\eta(\nabla_ZV)\}
  \end{eqnarray}
\medskip
\medskip
Now, letting $Y=Z=\xi$ gives rise to $(\alpha^2-\beta^2)\eta(\nabla_\xi V)=0$. Now, there will arise two cases either $\eta(\nabla_\xi V)=0$ or $(\alpha^2-\beta^2)=0$. From the definition of Ricci soliton (\ref{1.2}) we have,
\begin{equation}\label{3.11}
  \frac{1}{2}(g(\nabla_XV,Y)+g(\nabla_YV,X))+S(X,Y)=\lambda g(X,Y)
\end{equation}
for any vector fields $X$ and $Y$.For first case $\eta(\nabla_\xi V)=0$ which implies $\nabla_\xi V$ is orthogonal to $\xi$, putting $X=Y=\xi$ in (\ref{3.11}) gives $2(\alpha^2-\beta^2)=\lambda$. It directly implies that the soliton is shrinking if $\alpha^2<\beta^2$, steady if $\alpha^2=\beta^2$ and expanding if $\alpha^2>\beta^2$.\par
 For the second case where $\alpha^2=\beta^2$, then it follows directly from (\ref{3.11}) that $\nabla_\xi V=\lambda\xi$ i.e., the covariant derivative of the potential vector field $V$ in the direction of $\xi$ is a $\lambda$ multiple of $\xi$.
\end{proof}
\end{thm}

\begin{thm}
Let the $M$ be a 3-dimensional trans-Sasakian manifold of type $(\alpha,\beta)$ where the structure functions $\alpha$ and $\beta$ are constant with $\alpha\neq0$. If the metric $g$ represents a $*$-conformal Ricci soliton then the scalar curvature of the manifold is given by $r=(1-\frac{\beta^2}{\alpha^2})(\frac{p}{2}+\frac{1}{3}-\lambda+4\alpha^2)$.
\begin{proof}
  Since the metric $g$ represents a $*$-conformal Ricci solition, using (\ref{2.18}) in the definition of $*$-conformal Ricci solition (\ref{1.4}) we get,
  \begin{equation}\label{4.1}
    (\mathcal{L}_Vg)(X,Y)=(p+\frac{2}{3}+4(\alpha^2-\beta^2)-r-2\lambda)g(X,Y)+(r-4(\alpha^2-\beta^2))\eta(X)\eta(Y)
  \end{equation}
  for all vector fields $X$ and $Y$ on $M$. If we consider covariant derivative w.r.t. arbitrary vector field $Z$ then (\ref{4.1}) reduces to,
  \begin{eqnarray}\label{4.2}
    (\nabla_Z\mathcal{L}_Vg)(X,Y)&=&(Zr)[\eta(X)\eta(Y)-g(X,Y)]-(r-4(\alpha^2-\beta^2))\nonumber\\
    &&[\alpha g(\phi Z,X)\eta(Y)-\beta g(\phi X,\phi Z)\eta(Y)+\alpha g(\phi Z,Y)\nonumber\\
    &&\eta(X)-\beta g(\phi Y,\phi Z)\eta(X)]
  \end{eqnarray}
  for all $X,Y,Z\in\chi(M)$. Again from Yano\cite{Yano} we have the following commutation formula,
  \begin{eqnarray*}
    (\mathcal{L}_V\nabla_Zg-\nabla_Z\mathcal{L}_Vg-\nabla_{[V,Z]}g)(X,Y)&=&-g((\mathcal{L}_V\nabla)(X,Z),Y)\\
    &&-g((\mathcal{L}_V\nabla)(Y,Z),X),
  \end{eqnarray*}
   where $g$ is the metric connection i.e., $\nabla g=0$. So the above equation reduces to,
  \begin{equation}\label{4.3}
    (\nabla_Z\mathcal{L}_Vg)(X,Y)=g((\mathcal{L}_V\nabla)(X,Z),Y)+g((\mathcal{L}_V\nabla)(Y,Z),X),
  \end{equation}
   for all vector fields $X$, $Y$, $Z$ on $M$. Combining (\ref{4.2}) and (\ref{4.3}) and by a straightforward combinatorial computation and using the symmetry of $(\mathcal{L}_V\nabla)$ along with (\ref{2.7}) the foregoing equation yields,
   \begin{eqnarray}\label{4.4}
     (\mathcal{L}_V\nabla)(X,Y)&=&\frac{1}{2}(Dr)[g(X,Y)-\eta(X)\eta(Y)]-\frac{1}{2}(Xr)[Y-\eta(Y)\xi]-\nonumber\\
     &&\frac{1}{2}(Yr)[X-\eta(X)\xi]+(r-4(\alpha^2-\beta^2))[\beta g(\phi X,\phi Y)\xi\nonumber\\
     &&-\alpha\eta(Y)(\phi X)-\alpha\eta(X)(\phi Y)]
   \end{eqnarray}
   for arbitrary vector fields $X$ and $Y$ on $M$. Setting $Y=\xi$ in (\ref{4.4}) we have,
   \begin{equation}\label{4.5}
     (\mathcal{L}_V\nabla)(X,\xi)=-\frac{1}{2}(\xi r)[X-\eta(X)\xi]-\alpha(r-4(\alpha^2-\beta^2))(\phi X).
   \end{equation}
   Applying covariant derivative w.r.t. arbitrary vector field $Y$ and making use of (\ref{2.8}), (\ref{2.9}) and (\ref{2.10}) we obtain,
   \begin{eqnarray}\label{4.6}
     (\nabla_Y\mathcal{L}_V\nabla)(X,\xi) &=& \alpha(\mathcal{L}_V\nabla)(X,\phi Y)-\beta(\mathcal{L}_V\nabla)(X,Y)-\frac{1}{2}(Y(\xi r)) \nonumber\\
     && [X-\eta(X)\xi]+\frac{1}{2}(\xi r)[\alpha g(\phi X,Y)\xi+\beta g(\phi X,\phi Y)\xi\nonumber\\
     && -\alpha\eta(X)(\phi Y)+\beta\eta(X)Y-\beta\eta(Y)X]-\alpha(Yr)(\phi X)\nonumber\\
     &&-\alpha(r-4(\alpha^2-\beta^2))[\alpha g(X,Y)\xi-\alpha\eta(X)Y+\nonumber\\
     &&\beta g(\phi Y,X)\xi-\beta\eta(X)(\phi Y)+\beta\eta(Y)(\phi X)].
   \end{eqnarray}
   From Yano (for more details see \cite{Yano}) we know $(\mathcal{L}_VR)(X,Y)Z=(\nabla_X\mathcal{L}_V\nabla)\\(Y,Z)-(\nabla_Y\mathcal{L}_V\nabla)(X,Z)$. Using (\ref{4.6}) in this formula we get,
   \begin{eqnarray*}
     (\mathcal{L}_VR)(X,Y)\xi &=& \alpha(\mathcal{L}_V\nabla)(\phi X,Y)-\alpha(\mathcal{L}_V\nabla)(X,\phi Y)-\frac{1}{2}(X(\xi r))[Y- \\
     && \eta(Y)\xi]+\frac{1}{2}(Y(\xi r))[X-\eta(X)\xi]+\frac{1}{2}(\xi r)[2\alpha g(X,\phi Y)\xi \\
     && -\alpha\eta(Y)(\phi X)+\alpha\eta(X)(\phi Y)+2\beta\eta(Y)X-2\beta\eta(X)Y]\\
     && -\alpha(Xr)(\phi Y)+\alpha(Yr)(\phi X)-\alpha(r-4(\alpha^2-\beta^2))\\
     &&[\alpha\eta(X)Y-\alpha\eta(Y)X+2\beta g(\phi X,Y)\xi+2\beta\eta(X)(\phi Y)\\
     &&-2\beta\eta(Y)(\phi X)]
   \end{eqnarray*}
   Setting $Y=\xi$ in the foregoing equation we get,
   \begin{eqnarray}\label{4.7}
     (\mathcal{L}_VR)(X,\xi)\xi &=& \frac{1}{2}(\xi(\xi r))[X-\eta(X)\xi]+\beta(\xi r)[X-\eta(X)\xi]-\nonumber\\
     &&2\alpha(r-4(\alpha^2-\beta^2))[-\alpha X+\alpha\eta(X)\xi-\beta(\phi X)\nonumber\\
     &&
   \end{eqnarray}
   Again, Lie differentiation of the equation (\ref{2.14}) w.r.t. soliton vector field $V$ and use of (\ref{2.11}) and (\ref{2.14}) leads to,
   \begin{equation}\label{4.8}
     (\mathcal{L}_VR)(X,\xi)\xi=(\alpha^2-\beta^2)[g(X,\mathcal{L}_V\xi)\xi-((\mathcal{L}_V\eta)X)\xi-2\eta(\mathcal{L}_V\xi)X]
   \end{equation}
   which holds good for arbitrary vector field $X$ on $M$. Setting $Y=\xi$ in (\ref{4.1})implies,
   \begin{equation}\label{4.9}
     (\mathcal{L}_V\eta)X-g(X,\mathcal{L}_V\xi)=(p+\frac{2}{3}-2\lambda)\eta(X).
   \end{equation}
   Taking (\ref{4.9}) into account, Lie derivative of $\eta(\xi)=1$ along the direction of $V$ leads to
   \begin{equation}\label{4.10}
     2\eta(\mathcal{L}_V\xi)=-(p+\frac{2}{3}-2\lambda).
   \end{equation}
   After using (\ref{4.9}) and (\ref{4.10}), the equation (\ref{4.8}) reduces to,
   \begin{equation}\label{4.11}
     (\mathcal{L}_VR)(X,\xi)\xi=(\alpha^2-\beta^2)(p+\frac{2}{3}-2\lambda)[X-\eta(X)\xi]
   \end{equation}
   for all $X\in\chi(M)$
   Comparing (\ref{4.7}) with (\ref{4.11}) leads to,
   \begin{eqnarray}\label{4.12}
     &&(\alpha^2-\beta^2)(p+\frac{2}{3}-2\lambda)[X-\eta(X)\xi]=\frac{1}{2}(\xi(\xi r))[X-\eta(X)\xi]+\beta(\xi r)\nonumber\\
     &&[X-\eta(X)\xi]-2\alpha(r-4(\alpha^2-\beta^2))[-\alpha X+\alpha\eta(X)\xi-\beta(\phi X)
   \end{eqnarray}
   for $X\in\chi(M)$. Inner product w.r.t. arbitrary vector field $Y$ gives,
   \begin{eqnarray}\label{4.13}
     &&[\frac{1}{2}(\xi(\xi r))+\beta(\xi r)+2\alpha^2(r-4(\alpha^2-\beta^2))-(\alpha^2-\beta^2)(p+\frac{2}{3}-2\lambda)]\nonumber\\
     &&[g(X,Y)-\eta(X)\eta(Y)]+2\alpha\beta(r-4(\alpha^2-\beta^2))g(\phi X,Y)=0.
   \end{eqnarray}
   Anti-symmetrizing the foregoing equation yields,
   \begin{equation}\label{4.14}
     [\frac{1}{2}(\xi(\xi r))+\beta(\xi r)+2\alpha^2(r-4(\alpha^2-\beta^2))-(\alpha^2-\beta^2)(p+\frac{2}{3}-2\lambda)]g(\phi X,\phi Y)=0.
   \end{equation}
   Since this equation holds for arbitrary vector fields $\phi X$ and $\phi Y$ and as we know from (\ref{3.3}) that $\xi r=-2r\beta+12\beta(\alpha^2-\beta^2)$ holds in 3-dimensional trans-Sasakian manifold, we can easily conclude that the scalar curvature of the manifold satisfies $r=(1-\frac{\beta^2}{\alpha^2})(\frac{p}{2}+\frac{1}{3}-\lambda+4\alpha^2)$.
\end{proof}
\end{thm}

\section{\textbf{Example of 3-dimensional trans-Sasakian manifold admitting Ricci soliton }}
In this section we provide an example to verify our outcomes.
\begin{ex}
       We consider the manifold as $M=\{(x,y,z)\in\mathbb{R}^3:y\neq 0\}$, where $(x,y,z)$ are the standard coordinates in $\mathbb{R}^3$. The vector fields as defined bellow:
\begin{align*}
  e_1 &=e^{2z}\frac{\partial}{\partial x}, & e_2&=e^{2z}\frac{\partial}{\partial y}, & e_3=\frac{\partial}{\partial z}
\end{align*}
are linearly independent at each point of $M$. The Riemannian metric $g$ is defined by:
\[
  g_{ij}=\begin{cases}
               1, & \mbox{if } i=j~and~i,j\in\{1,2,3\} \\
               0, & \mbox{otherwise.}
             \end{cases}
  \]
Let $\xi=e_3$. Then the 1-form $\eta$ is defined by $\eta(Z)=g(Z,e_3)$, for arbitrary $Z\in\chi(M)$, then we have the following relations:
\begin{align*}
  \eta(e_1)=\eta(e_2) &= 0, & \eta(e_3) &= 1.
\end{align*}
Let us define the (1,1)-tensor field $\phi$ as
\begin{align*}
  \phi e_1 &=e_2, & \phi e_2 &= -e_1, & \phi e_3 &= 0,
\end{align*}
 then it satisfies,
\begin{eqnarray}
  \nonumber \phi^2(Z)&=&-Z+\eta(Z)e_3,\\
  \nonumber g(\phi Z,\phi W)&=&g(Z,W)-\eta(Z)\eta(W)
\end{eqnarray}
for arbitrary $Z,W\in\chi(M)$.

Thus $(\phi,\xi,\eta,g)$ defines an almost contact metric structure on $M$. We can now easily conclude:
\begin{align*}
  [e_1,e_2] &= 0, & [e_2,e_3] &= -2e_2, & [e_1,e_3] &=-2e_1.
\end{align*}
Let $\nabla$ be the Levi-Civita connection of $g$. Then from $Koszul's formula$ for arbitrary $X,Y,Z\in\chi(M)$ given by:\\
\begin{eqnarray}
2g(\nabla_XY,Z)&=&Xg(Y,Z)+Yg(Z,X)-Zg(X,Y)-g(X,[Y,Z])\nonumber\\
&-&g(Y,[X,Z])+g(Z,[X,Y]),\nonumber
\end{eqnarray}
we can have:
\begin{align*}
  \nabla_{e_1}e_1 &=2e_3, & \nabla_{e_1}e_2 &= 0, & \nabla_{e_1}e_3 &=-2e_1, \\
  \nabla_{e_2}e_1 &=0 ,& \nabla_{e_2}e_2 &=2e_3, & \nabla_{e_2}e_3 &=-2e_2, \\
  \nabla_{e_3}e_1 &=0, & \nabla_{e_3}e_2 &=0 ,& \nabla_{e_3}e_3 &=0.
\end{align*}
From here we can easily verify that the relations (\ref{2.8}) and (\ref{2.9}) are satisfied. Hence the considered manifold is trans-Sasakian manifold of type $(0,-2)$. The components of Riemannian curvature tensor are given by,
\begin{align*}
  R(e_1,e_2)e_1&=-4e_3, & R(e_1,e_2)e_2&=-4e_1, & R(e_1,e_2)e_3&=0, \\
  R(e_1,e_3)e_1&=4e_2, & R(e_1,e_3)e_2&=0 ,& R(e_1,e_3)e_3&=-4e_1, \\
  R(e_2,e_3)e_1&=0, & R(e_2,e_3)e_2&=-4e_2, & R(e_2,e_3)e_3&=-4e_2.
\end{align*}
And the components of Ricci tensor and $*$-Ricci tensor are given by:
\begin{align*}
  S(e_1,e_1)&=0, & S(e_2,e_2)&=0, & S(e_3,e_3)&=-8.
\end{align*}
From here we can easily deduce that the scalar curvature of the manifold $r=\sum_{i=1}^{3}S(e_i,e_i)=-8$.
Let us define a vector field by, $V=\xi$. Then we can obtain:
\begin{align*}
  (\mathcal{L}_Vg)(e_1,e_1)&=-4, & (\mathcal{L}_Vg)(e_2,e_2)&=-4, & (\mathcal{L}_Vg)(e_3,e_3)&=0.
\end{align*}
Contracting (\ref{1.2}) and using the result $r=-8$ we deduce $\lambda=4$. So $g$ defines a Ricci soliton on this trans-Sasakian manifold for $\lambda=4$.
\end{ex}

\end{document}